\documentclass[10pt,twoside]{article}
 \usepackage{amssymb,latexsym,amsmath,amsthm}
 
 \usepackage{fullpage}

\newtheorem{thm}{Theorem}
\newtheorem{lemma}{Lemma}

\newtheorem{cor}{Corollary}

\newtheorem{rem*}{Remark}

\newtheorem{conj}{Conjecture}

\begin{document}
\date{}

\title{On the H\'{e}non-Lane-Emden conjecture}
\author{Mostafa Fazly\thanks{Research partially supported by a University Graduate Fellowship at the University of British Columbia.} \quad and \quad Nassif  Ghoussoub\thanks{Partially supported by a grant
from the Natural Sciences and Engineering Research Council of Canada.  }
\quad   
\\
\small Department of Mathematics,\\
\small University of British Columbia, \\
\small Vancouver BC Canada V6T 1Z2 \\
\small {\tt fazly@math.ubc.ca}\\
\small {\tt nassif@math.ubc.ca} \\\\
\today\\
}
\maketitle

\vspace{3mm}

\begin{abstract}
We consider Liouville-type theorems for the following H\'{e}non-Lane-Emden system 
\begin{eqnarray*}
 \left\{ \begin{array}{lcl}
\hfill -\Delta u&=& |x|^{a}v^p   \ \ \text{in}\ \ \mathbb{R}^N,\\   
\hfill -\Delta v&=& |x|^{b}u^q   \ \ \text{in}\ \ \mathbb{R}^N,
\end{array}\right.
  \end{eqnarray*}
when $pq>1$, $p,q,a,b\ge0$. The main conjecture states that there is no non-trivial non-negative solution whenever  $(p,q)$ is  {\it under} the critical Sobolev hyperbola, i.e.
$ \frac{N+a}{p+1}+\frac{N+b}{q+1}>{N-2}$.
 We show that this is indeed the case  in dimension $N=3$ provided the solution is also assumed to be bounded, extending a result established recently by Phan-Souplet  in the scalar case.
 
Assuming stability of the solutions, we could then prove Liouville-type theorems in higher dimensions. 
  For the scalar cases, albeit of second order ($a=b$ and $p=q$) or of fourth order ($a\ge 0=b$ and $p>1=q$), we show that for all dimensions $N\ge 3$ in the first case (resp., $N\ge 5$ in the second case), there is no positive solution with a finite Morse index, whenever $p$ is below the corresponding critical exponent, i.e $  1<p<\frac{N+2+2a}{N-2}$ (resp., $  1<p<\frac{N+4+2a}{N-4}$). Finally, we show that non-negative stable solutions of the full H\'{e}non-Lane-Emden system are trivial provided  
  \begin{equation*}\label{sysdim00}
 N<2+2\left(\frac{p(b+2)+a+2}{pq-1}\right) \left(\sqrt{\frac{pq(q+1)}{p+1} }+ \sqrt{ \frac{pq(q+1)}{p+1}-\sqrt\frac{pq(q+1)}{p+1}  } \right).
 \end{equation*}

\end{abstract}


\section {Introduction and main results}
We consider the following weighted system
\begin{eqnarray}
\label{mainbound}
 \left\{ \begin{array}{lcl}
\hfill -\Delta u&=& |x|^{a}v^p   \ \ \text{in}\ \ \mathbb{R}^N,\\   
\hfill -\Delta v&=& |x|^{b}u^q   \ \ \text{in}\ \ \mathbb{R}^N,
\end{array}\right.
  \end{eqnarray}
where $pq>1$ and $p,q,a,b\ge0$ and $\Omega$ is a subset of $\mathbb{R}^N$, $N\ge 1$. 

We start by noting that in the case of  the Lane-Emden scalar equation (i.e., when $p=q$ and $a=b=0$) on a bounded star-shaped domain $\Omega\subset\mathbb{R}^N$, the Pohozaev inequality shows that there is no positive solution satisfying the Dirichlet boundary condition, whenever  $p \geq \frac{N+2}{N-2}$, the critical Sobolev exponent.  On the other hand, a celebrated theorem by Gidas-Spruck  \cite{gs}  states that there is no positive solution for the Lane-Emden equation on the whole space whenever $p< \frac{N+2}{N-2}$ for $N\ge 3$.   This non-existence result is also optimal as shown by Gidas, Ni and Nirenberg in \cite{gnn} under the assumption that $u=O(|x|^{2-N})$, and by Caffarelli, Gidas and Spruck in \cite{cgs} without  the growth assumption. See also Chen and Li \cite{cl}  for an easier proof based on the  moving planes method. Also, Lin \cite{lin} using moving plane methods proved similar optimal non-existence results for $p< \frac{N+4}{N-4}, N>4$ in the case of  the fourth order Lane-Emden equation (i.e., when $p>1=q$ and $a=b=0$).

In the case of the system (\ref{mainbound}), one can again use the Pohozaev identity whenever $\Omega$ is a bounded star-shaped domain in $\mathbb{R}^N$, to  establish the following non-existence result.

\medskip
\noindent{\bf Theorem A. \cite{mostafa,pserrin}}\quad
{\it Let  $N\ge3$ and let $\Omega\subset\mathbb{R}^N$ be a star-shaped bounded domain. If 
 \begin{equation}\label{curvebound}
   \frac{N+a}{p+1}+\frac{N+b}{q+1}\le{N-2},
  \end{equation}
}
then there is no positive solution for (\ref{mainbound}) on $\Omega$ that satisfy the Dirichlet boundary conditions.
\medskip

By noting that the curve $ \frac{N+a}{p+1}+\frac{N+b}{q+1}={N-2}$ is the {\it critical Sobolev hyperbola}, the above theorem states that the Liouville-type result for positive solutions on bounded star-shaped domain holds when $(p, q)$ is {\it above} the critical hyperbola.  It  is therefore expected that -- just like the case of  the scalar Lane-Emden equation ($p=q$ and $a=b=0$) -- the non-existence of solutions on the whole space $\mathbb{R}^N$ should occur exactly when $(p,q)$ is in the complementary domain, that is when it is {\it under} the critical hyperbola.

This is the statement of the following  {\it H\'{e}non-Lane-Emden conjecture}.

\begin{conj} Suppose $(p,q)$ is {\it under} the critical hyperbola, i.e.,
  \begin{equation}\label{curve}
   \frac{N+a}{p+1}+\frac{N+b}{q+1}>{N-2}.
  \end{equation}
Then there is no positive solution for system (\ref{mainbound}).
\end{conj}

Proving such a non-existence result seems to be challenging even for the Lane-Emden conjecture (i.e., when $a=b=0$) for systems.  The case of radial solutions was solved by Mitidieri \cite{m} in any dimension, and both Mitidieri \cite{m} and Serrin-Zou \cite{sz98} constructed positive radial solutions {\it on} and {\it above} the critical hyperbola, i.e. $\frac{1}{p+1}+\frac{1}{q+1}\le\frac{N-2}{N}$, which means that the non-existence theorem is optimal for radial solutions.\\
 For non-radial solutions of the Lane-Emden system, there are  the results of Souto \cite{s95}, Mitidieri \cite{m} and Serrin-Zou \cite{sz} who proved the non-existence of solutions in dimensions $N=1,2$, while in  dimension $N=3$,  Serrin-Zou \cite{ sz} gave a proof for the non-existence of polynomially bounded solutions, an assumption that was removed later by Pol\'{a}\v{c}ik, Quittner and Souplet \cite{pqs}. More recently, Souplet \cite{so} settled completely the conjecture in dimension $N=4$, while providing in dimensions $N\ge 5$, a more restrictive new region for the exponents $(p,q)$ that insures non-existence. 
 
  \medskip
\noindent{\bf Theorem B. (Souplet \cite{so})}\quad Assume $a=b=0$.
{\it  
\begin {enumerate}
\item[(i)] Let $N=4$ and $p, q >0$. If $(p, q)$ satisfies 
\begin{equation}\label{curvesouplet}
\frac{1}{p+1}+\frac{1}{q+1}>\frac{N-2}{N},
\end{equation}
then system (\ref{mainbound}) has no positive  solutions.

\item [(ii)] Let $N\ge 5$, and $p, q >0$ with $pq > 1$. If  $(p, q)$ satisfies (\ref{curvesouplet}), along with
\begin{equation}
2 \max\left\{\frac{p+1}{pq-1},\frac{q+1}{pq-1}\right\} > N-3,
\end{equation}
then every non-negative solution of system (\ref{mainbound}) is necessarily trivial.

\end{enumerate}

}

The Lane-Emden conjecture in dimensions $N\ge 5$ is still open. The H\'{e}non-Lane-Emden conjecture is even less understood. Even for the scalar case $a=b$ and $p=q$ ( i.e., the H\'{e}non equation), Gidas and Spruck in \cite{gs} solved the conjecture only for radial solutions, also showing that in this case, the non-existence result is optimal. For non-radial solutions, they proved some partial results such as the non-existence of positive solutions for $a\ge 2$ and $p\le \frac{N+2}{N-2}$ (the Sobolev critical exponent for $a=0$). 

For systems,  Mitidieri \cite{m} gave a partial solution to the conjecture for radial solutions by showing that there is no positive radial solution for (\ref{mainbound}) for all $N\ge3$ provided $p,q>1$ satisfy
 \begin{equation}\label{radialcurvemin}
   \frac{N+\min\{a,b\}}{p+1}+\frac{N+\min\{a,b\}}{q+1}>{N-2},
  \end{equation}
Recently, Bidaut-Veron-Giacomini \cite{bvg} used a Pohozaev type argument and a suitable change of variables to  give a complete solution in the radial case. 

\medskip
\noindent{\bf Theorem C. (Bidaut-Veron-Giacomini  \cite{bvg})}\quad
{\it For  $N\ge3$, System  (\ref{mainbound}) admits a positive radial solution $(u,v)$ such that $u,v\in C^{2}(0,\infty)\cap C([0,\infty))$ if and only if $(p,q)$ is above or on the critical hyperbola, i.e., when (\ref{curvebound}) holds.
}
\medskip

\subsection{Liouville theorems for bounded non-negative solutions}

With the lack of progress on the full conjecture, the attention turned to showing that bounded non-negative solutions are necessarily trivial.  Recently, Phan and Souplet \cite{ps} showed among other things that  the H\'{e}non-Lane-Emden conjecture for the scalar case holds for bounded positive solutions in dimension $N=3$. 

  \medskip
\noindent{\bf Theorem D. (Phan-Souplet \cite{ps})}\quad
{\it Let $N=3$, $a=b>0$ and $p=q>1$. Assume $(p,q)$ satisfies (\ref{curve}), then there is no positive bounded solution for the H\'{e}non equation, i.e., 
\begin{equation}\label{mainsingle}
  -\Delta u=|x|^a u^p\, \ \ \text{in} \ \ \mathbb{R}^N.
  \end{equation} 
}
In this note, we shall first extend the above result of  Phan-Souplet \cite{ps} to the full H\'{e}non-Lane-Emden system by showing the following\footnote{Upon receiving our preprint, P. Souplet informed us that Q.H. Phan has also proved the same result in dimension $N=3$, as well as other interesting results in higher dimensions. Our proofs are quite similar since both are essentially refinements of those of P. Souplet in his groundbreaking work on the Lane-Emden conjecture for systems.}. 
\begin{thm}\label{result} Suppose $N=3$ and $(p,q)$ satisfy (\ref{curve}). Then, 
there is no positive bounded solution for (\ref{mainbound}).
\end{thm}

 We also give a few  partial results for the H\'{e}non equation whether of second order or fourth order in all dimensions $N\ge 3$ or $N\ge 5$.

We note that  Miditieri and Pohozaev \cite{mp} have shown that the above result holds in higher dimension provided the following stronger condition holds:
\begin{equation*}
   \max\{\alpha,\beta\}\ge N-2,
  \end{equation*}
where $\alpha:=\frac{(b+2)p+(a+2)}{pq-1}$ and $\beta:=\frac{(a+2)q+(b+2)}{pq-1}$. For that they used a rescaled test-function method (as in Lemma \ref{weightedest} below)  to prove the result for $p,q\ge 1$. More recently, Armstrong and Sirakov \cite{as} proved --among other things-- similar results for $p,q>0$, by developing new maximum principle type arguments. We are thankful to P. Souplet for informing us of these latest developments by Armstrong and Sirakov.

\subsection{Liouville theorems for stable non-negative solutions}

We shall also consider in the scalar case the question of existence of solutions with finite Morse index solutions (as opposed to bounded solutions). For scalar equations, we get the following counterpart to the Phan-Souplet result in  higher dimensions ($N\ge 3$).

\begin{thm}\label{resultMorse}
Let $a\ge0$, $p>1$ and $N\ge3$. Then, for any Sobolev sub-critical exponent, i.e.,
$$  1<p<\frac{N+2+2a}{N-2}, $$
equation (\ref{mainsingle}) has no positive solution with finite Morse index.
\end{thm}

We also have the following result for the fourth order equation,
\begin{equation}\label{fourth}
\Delta^2 u=|x|^a u^p \ \ \ \ \text{in} \ \ \mathbb{R}^N.
\end{equation}

\begin{thm}\label{fourthresultMorse}
Let $a\ge0$, $p>1$ and $N\geq 5$. Then, for any Sobolev sub-critical exponent, i.e.,
$$  1<p<\frac{N+4+2a}{N-4}, $$
equation (\ref{fourth}) has no positive solution with finite Morse index.
\end{thm}

For systems, we have the following result.

\begin{thm}\label{system}
Suppose that $0\le a-b\le (N-2)(p-q)$. Then, system (\ref{mainbound})  has no positive stable solution whenever the dimension satisfy
\begin{equation}\label{sysdim}
 N<2+2\left(\frac{p(b+2)+a+2}{pq-1}\right) \left(\sqrt{\frac{pq(q+1)}{p+1} }+ \sqrt{ \frac{pq(q+1)}{p+1}-\sqrt\frac{pq(q+1)}{p+1}  } \right).
 \end{equation}
\end{thm}
The case when $a=b=0$ (i.e., the Lane-Emden system) was already established by by Cowan in \cite{cowan}.
Note that this result contains the result of Fazly in  \cite{mostafa}, who had considered the case $q =1<p $,  $a=b$ and shown the result under the condition,
\begin{equation}\label{sysfazly}
N<8+3a+\frac{8+4a}{p-1},
\end{equation}
which is already larger than the domain under the critical hyperbola, i.e. $N< 4+a+\frac{8+4a}{p-1}$.  Also, this  contains  the result of Wei-Ye in \cite{wy} who had considered the case $q =1<p $,  $a=b=0$. There are also various results for the cases where  $-2<a,b<0$ and $pq\le 1$. For that we refer to  \cite{bvg,mostafa,mp,gy,gs,gs2,ps}.

\section{Proof in the case of non-negative solutions}

In this section, we shall prove here Theorem \ref{result}. The main tools will be Pohozaev-type identities as well as various  integral estimates.

The proof is heavily inspired by ideas of Souplet \cite{so} and Serrin-Zou \cite{sz}. We use Pohozaev-type identities, various integral estimates, as well as some elliptic estimates on the sphere. Throughout this section, all norms refer to functions defined  on the unit sphere, i.e. $||u||_{m}:=||u||_{L^{m}(S^{N-1})}$.

We start with the following estimate on the non-linear terms. Note that for $a=b=0$, this was proved by Serrin and Zou \cite{sz} via ODE techniques, and by Miditieri
and Pohozaev \cite{mp} who used the following rescaled test functions approach for $a,b>-2$. For the sake of convenience of readers, we recall the proof. Interested readers can find more details for both scalar and system cases in \cite{qs}.
\begin{lemma} \label{weightedest}
For any positive entire solution $(u,v)$ of (\ref{mainbound}) and $R>1$, there holds  
\begin{eqnarray} \label{weightedest1}
\int_{B_R}{       |x|^{a}   v^p} &\le&  C \ R^{N -2- \frac{(b+2)p+(a+2)}{pq-1} },
\\ \label{weightedest2}
\int_{B_R}{    |x|^{b}   u^q}  &\le& C\  R^{N-2 - \frac{(a+2)q+(b+2)}{pq-1} },
 \end{eqnarray}
where the positive constant $C$ does not depend on $R$.
\end{lemma}   

\textbf{Proof:} Fix the following function $\zeta_R\in C^2_c(\mathbb{R}^N)$ with $0\le\zeta_R\le1$;
 $$\zeta_R(x)=\left\{
                      \begin{array}{ll}
                        1, & \hbox{if $|x|<R$;} \\
                        0, & \hbox{if $|x|>2R$;} 
                                                                       \end{array}
                    \right.$$
where $||\nabla\zeta_R||_{\infty} \le \frac{C}{R}$ and $||\Delta\zeta_R||_\infty\le  \frac{C}{R^2}$. For fixed $m\ge 2$, we have  

$$|\Delta \zeta^m_R(x)|\le C \left\{
                      \begin{array}{ll}
                        0, & \hbox{if $|x|<R$ or $|x|>2R$;} \\
                         R^{-2} \zeta^{m-2}_R , & \hbox{if $R<|x|<2R$;} 
                                                                       \end{array}
                    \right.$$
For $m\ge 2$, test the first equation of (\ref{mainbound}) by $\zeta^m_R$ and integrate to get 
\begin{eqnarray*}
 \int_{\mathbb{R}^N}     |x|^{a}     v^p  \zeta^m_R &= &- \int_{\mathbb{R}^N}  \Delta u   \zeta^m_R\\
&=& - \int_{\mathbb{R}^N} u\Delta \zeta^m_R   \le C  R^{-2}\int_{B_{2R}\setminus B_{R}} u\zeta^{m-2}_R.
\end{eqnarray*}
Applying H\"{o}lder's inequality we get 
\begin{eqnarray*}
\int_{\mathbb{R}^N}      |x|^{a}     v^p  \zeta^m_R & \le  & C \ R^{-2}  \left(         \int_{B_{2R}\setminus B_{R}}  |x|^{\frac{-b}{q }q'    }          \right)^{\frac{1}{q'}}   \left(     \int_{B_{2R}\setminus B_{R}}      |x|^{b}    u^q \zeta^{(m-2)q}_R   \right)^{1/q}\\
& \le  & C\  R^{  (N-\frac{b}{q}q')\frac{1}{q'}  -2 }  \left(     \int_{B_{2R}\setminus B_{R}}       |x|^{b}     u^q \zeta^{(m-2)q}_R   \right)^{1/q}.
\end{eqnarray*}
By a similar calculation for $k\ge 2$, we obtain
\begin{eqnarray*}
 \int_{\mathbb{R}^N}    |x|^{b}      u^q  \zeta^k_R & \le &C \ R^{  (N-\frac{a}{p}p')\frac{1}{p'}  -2 }     \left(\int_{B_{2R}\setminus B_{R}}    |x|^{a} v^p \zeta^{(k-2)p}_R   \right)^\frac{1}{p} ,
 \end{eqnarray*}
where $\frac{1}{p}+\frac{1}{p'}=1$. Since $pq>1$, for large enough $k$ we have $2+\frac{k}{q}<(k-2)p$. So, we can choose $m$ such that $2+\frac{k}{q}\le   m  \le   (k-2)p$ which means that $m\le (k-2)p$ and $k\le (m-2)q$. By collecting  the above inequalities we get for $pq> 1$,
 
 \begin{eqnarray}
\nonumber   \left(\int_{\mathbb{R}^N}{         |x|^{a}      v^p   \zeta^m_R} \right)^{pq}&\le & C\    R^{ [ (N-\frac{b}{q}q')\frac{1}{q'}  -2 ] pq}  \left( \int_{B_R}{    |x|^{b}  u^q \zeta^k_R}  \right)^p 
\\ \label{L1final1}
&\le & C\ R^{   (N-2)(pq-1) -[(b+2)p+(a+2)] }  \int_{B_{2R}\setminus B_{R}}     |x|^{a}      v^p   \zeta^m_R,
 \end{eqnarray}
 and
 \begin{eqnarray}
  \nonumber       \left(\int_{\mathbb{R}^N}{    |x|^{b}  u^q     \zeta^k_R } \right)^{pq} &\le & C\    R^{ [  (N-\frac{a}{p}p')\frac{1}{p'}  -2  ] pq}   \left(\int_{B_{R}}{         |x|^{a}      v^p   \zeta^m_R} \right)^q
\\ \label{L1final2}
&\le & C\ R^{   (N-2)(pq-1) -[(a+2)q+(b+2)] }  \int_{B_{2R}\setminus B_{R}}      |x|^{b}  u^q     \zeta^k_R.
 \end{eqnarray}

\hfill $\Box$

By using H\"{o}lder's inequality, we can now get the following  $L^{1}$-estimates.
\begin{cor} \label{L1est}
With the same assumptions as Lemma \ref{L1est}, we have
 \begin{eqnarray*}
\int_{B_R}    v^s &\le&  C  R^{N-\frac{(a+2)q+(b+2)}{pq-1} s},\\
\int_{B_R}    u^t &\le&  C  R^{N-\frac{(b+2)p+(a+2)}{pq-1} t},
 \end{eqnarray*}
 for any $0<t<q$ and $0<s<p$ where the positive constant $C$ does not depend on $R$.
\end{cor}

We now recall the following fundamental elliptic estimates.

\begin{lemma} (Sobolev inequalities on the sphere $S^{N-1}$) \label{sobolev} 
Let $N\ge2$,  integer $j\ge 1$ and $1<k<m\le\infty$. For $z\in W^{j,k}(S^{N-1})$, we have
     $$||z||_{L^m(S^{N-1})}\le C( || D_\theta^j  z||_{L^k(S^{N-1})}  +|| z  ||_{L^1(S^{N-1})} ),$$ where
     
       $$\left\{
                      \begin{array}{ll}
                        \frac{1}{k}-\frac{1}{m}=\frac{j}{N-1}, & \hbox{if $k<(N-1)/j$,} \\
                        m=\infty, & \hbox{if $k>(N-1)/j$,} 
                                                                       \end{array}
                    \right.$$
                     and $C=C(j,k,N)>0$.
\end{lemma}

\begin{lemma} \label{ellip} 
(Elliptic $L^p$-estimate on $B_R$).  Let $1<k<\infty$ and $R>0$. For $z\in W^{2,k}(B_{2R})$, we have
$$\int_{B_R }  |  D_{x}^2 z |^k\le C\left(\int_{B_{2R} }  |\Delta z|^k  +   R^{-2k}    \int _{B_{2R}  }    |z|^k  \right ),$$
where $C=C(k,N)>0$.
\end{lemma}

\begin{lemma}\label{interp} 
(An interpolation inequality on $B_R$).  Let  $R>0$. For $z\in W^{2,1}(B_{2R})$, we have
$$\int_{B_R}    |  D_{x} z |\le C \left( R  \int_{B_{2R}  }  |\Delta z|  +    R^{-1}    \int _{B_{2R}  }    |z|  \right) ,$$

where $C=C(N)>0$.
\end{lemma} 

By applying Lemma \ref{weightedest}, Corollary \ref{L1est} and Lemma \ref{interp}, we obtain the following estimates on the derivatives of $u$ and $v$.
\begin{lemma} \label{DL1est}
We have
\begin{eqnarray*}
\int_{B_R}{    |  D_{x}v|}  &\le& C\  R^{N-1-\frac{(a+2)q+(b+2)}{pq-1} },
\\
\int_{B_R}{     | D_{x} u|} &\le&  C \ R^{N-1-\frac{(b+2)p+(a+2)}{pq-1} },
 \end{eqnarray*}
where the positive constant $C$ does not depend on $R$.
\end{lemma}

\begin{lemma}\label{regularity}
($L^1$-regularity estimate on $B_R$) Let $N> 2$ and $1\le k<\frac{N}{N-2}$. For any $z\in L^1(B_{2R})$ we have 
$$ || z||_{L^k(B_R)}\le C \left( R^{2+N(\frac{1}{k}-1)} || \Delta z||_{L^1(B_{2R})} + R^{N (\frac{1}{k}-1) } ||z||_{L^1(B_{2R})}  \right),$$
where $C=C(k,N)>0$.
\end{lemma}

For $a=b=0$, the following Pohozaev identity has been obtained by Mitidieri \cite{mi}, Serrin and Zou \cite{sz}. It has also been used by Souplet in \cite{so}.

\begin{lemma}\label{Poho}
(Pohozaev identity). Suppose $\lambda,\gamma\in\mathbb{R}$ satisfy $\lambda+\gamma=N-2$.  If $(u,v)$ is a positive solution of (\ref{mainbound}), then it necessarily satisfy
\begin{eqnarray*}
\left(\frac{N+a}{p+1} -\lambda \right) \int_{B_{R}} |x|^{a} v^{p+1} &+& \left(\frac{N+b}{q+1} -\gamma \right) \int_{B_{R}} |x|^{b} u^{q+1} \\ 
&=& R^{N+a} \int_{S^{N-1}}  \frac{v^{p+1}}{p+1} +R^{N+b} \int_{S^{N-1}}\frac{u^{q+1}}{q+1} 
+R^{N} \int_{S^{N-1}} \left( u_{r}v_{r} -R^{-2} u_{\theta}v_{\theta}  \right) \\
&&+ R^{N-1} \int_{S^{N-1}}(\lambda u_{r}v+\gamma v_{r}u).
\end{eqnarray*}

\end{lemma}

Now, we are in the position to prove Theorem \ref{result}.
\\

\noindent \textbf{Proof of Theorem \ref{result}:} Since $(p,q)$ satisfy (\ref{curve}), then we can choose $\lambda$ and $\gamma$ such that $\frac{N+a}{p+1}>\lambda$ and $\frac{N+b}{q+1}>\gamma$. Now, for all $R>0$ define
 $$F(R):=\left(\frac{N+a}{p+1} -\lambda \right) \int_{B_{R}} |x|^{a} v^{p+1} + \left(\frac{N+b}{q+1} -\gamma \right) \int_{B_{R}} |x|^{b} u^{q+1}.$$ 
From Lemma \ref{Poho}, we have 
\begin{equation}\label{F}
F(R)\le C \left( G_{1}(R) +G_{2}(R)  \right),
\end{equation}
where 
$$G_{1}(R):=R^{N+a}\int_{S^{N-1}} v^{p+1}+R^{N+b}\int_{S^{N-1}} u^{q+1},$$ and $$G_{2}(R):=R^{N} \int_{S^{N-1}} 
\left(|D_{x}u(R)|+R^{-1}u(R)\right)\left(|D_{x}v(R)|+R^{-1}v(R)\right) .$$

\noindent \textbf{Step 1}. Upper bounds for $G_{1}$ and $G_{2}$. Set $m=\infty$ in Lemma \ref{sobolev} to get for either $t=p+1$ or $t=q+1$
$$|| u||_{t}\le|| u||_{\infty}\le C (|| D_{\theta}^{2} u||_{1+\epsilon}+||u||_{1})\le  C (R^{2}|| D_{x}^{2} u||_{1+\epsilon}+||u||_{1}),$$
where $\epsilon>0$ is small enough and will be chosen later. So,
\begin{eqnarray}
\nonumber   G_{1}(R) &\le& R^{N+a+2(p+1)} \left(  || D_{x}^{2} v||_{1+\epsilon} + R^{-2}|| v ||_{1}\right)^{1+p}\\ \label{G1}
&& +R^{N+b+2(q+1)} \left(  || D_{x}^{2} u||_{1+\epsilon} + R^{-2}|| u ||_{1}\right)^{1+q}.
\end{eqnarray}

We now look for the same type bounds for $G_{2} $. Apply Schwarz's inequality to get 
\begin{eqnarray*}
G_{2}(R) &\le& R^{N} \left(\int_{S^{N-1}} 
\left(|D_{x}u(R)|+R^{-1}u(R)\right)^{2} \right )^{1/2}      \left(\int_{S^{N-1}}   \left(|D_{x}v(R)|+R^{-1}v(R)\right)^{2} \right) ^{1/2}\\
&\le& R^{N} \left(||D_{x}u||_{2}+R^{-1}||u||_{1}\right)\left(||D_{x}v||_{2}+R^{-1}||v||_{1}\right).
\end{eqnarray*}
Then, using Lemma \ref{sobolev} we obtain the following upper bounds. 
\begin{eqnarray*}
||D_{x}u||_{2} &\le& C \left( || D_{\theta}D_{x}u||_{1+\epsilon} +||D_{x}u||_{1}  \right) \le C \left( R ||D_{x}^{2}u||_{1+\epsilon} +||D_{x}u||_{1}  \right), \\
||D_{x}v||_{2} &\le& C \left( || D_{\theta}D_{x}v||_{1+\epsilon} +||D_{x}v||_{1}  \right) \le C \left( R ||D_{x}^{2}v||_{1+\epsilon} +||D_{x}v||_{1}  \right). 
\end{eqnarray*}
It follows that
\begin{eqnarray}\label{G2}
G_{2}(R) &\le& R^{N+2} 
 \left(||D_{x}^{2}u||_{1+\epsilon}+R^{-1}||D_{x}u||_{1}+R^{-2}||u||_{1}\right)\left(||D_{x}^{2}v||_{1+\epsilon}+R^{-1}||D_{x}v||_{1}+R^{-2}||v||_{1}\right).
\end{eqnarray}

\noindent\textbf{Step 2}. The following $L^{t}$-estimates hold in the annulus domain $B_R\setminus B_{R/2}$;
\begin{eqnarray}
\label{v}\int_{R/2}^{R} || v(r)||_{1} r^{N-1} dr &\le& C\ R^{N-\frac{(a+2)q+(b+2)}{pq-1}},\\
\label{u}\int_{R/2}^{R} || u(r)||_{1} r^{N-1} dr  &\le& C\ R^{N-\frac{(b+2)p+(a+2)}{pq-1}},\\
\label{Dv}\int_{R/2}^{R} ||D_{x} v||_{1} r^{N-1} dr &\le& C\ R^{N-1-\frac{(a+2)q+(b+2)}{pq-1}},\\
\label{Du}\int_{R/2}^{R} || D_{x}u||_{1} r^{N-1} dr  &\le& C\ R^{N-1-\frac{(b+2)p+(a+2)}{pq-1}},\\
\label{D2v}\int_{R/2}^{R} ||D_{x}^{2} v||^{1+\epsilon}_{1+\epsilon} r^{N-1} dr &\le& C\ R^{N-2-\frac{(a+2)q+(b+2)}{pq-1}+b\epsilon},\\
\label{D2u}\int_{R/2}^{R} ||D_{x}^{2} u||^{1+\epsilon}_{1+\epsilon} r^{N-1} dr &\le& C\ R^{N-2-\frac{(b+2)p+(a+2)}{pq-1}+a\epsilon}.
\end{eqnarray}
To prove (\ref{v})-(\ref{Du}), we just apply Corollary \ref{L1est} and Lemma \ref{DL1est}. Here is for example the proof for (\ref{D2u}).  Apply Lemma \ref{ellip}, Corollary \ref{L1est} and Lemma \ref{weightedest} to get 
\begin{eqnarray*}
\int_{R/2}^{R}||D_{x}^{2} u||^{1+\epsilon}_{1+\epsilon} r^{N-1} dr & =& \int_{R/2}^{R} |D^{2}_{x}u|^{1+\epsilon}
dx \\
&\le& C \int _{B_{2R}} |\Delta u|^{1+\epsilon} dx +C\ R^{-2(1+\epsilon)}  \int _{B_{2R}} u^{1+\epsilon} dx\\
&\le& C\ R^{a\epsilon}\int _{B_{2R}} |x|^{a} v^{p(1+\epsilon)}dx+ C\ R^{-2(1+\epsilon)}\int _{B_{2R}} u\\
&\le& C\ R^{ N-2-\frac{(b+2)p+(a+2)}{pq-1}+a\epsilon }+ C\ R^{   N-\frac{(b+2)p+(a+2)}{pq-1} -2(1+\epsilon)}\\
&\le &C\ R^{ N-2-\frac{(b+2)p+(a+2)}{pq-1}+a\epsilon}.
\end{eqnarray*}
The proof of (\ref{D2v}) is similar.

\noindent\textbf{Step 3} For large enough $M$, define following sets;
\begin{eqnarray*}
\label{mv}\Gamma_{1}(R) &:=&\{r\ \in(R,2R); \  || v(r)||_{1}> M R^{-\frac{(a+2)q+(b+2)}{pq-1}}\},\\
\label{mu}\Gamma_{2}(R) &:=&\{r\ \in(R,2R); \  || u(r)||_{1}  > M R^{-\frac{(b+2)p+(a+2)}{pq-1} }\},\\
\label{mDv}\Gamma_{3}(R) &:=&\{r\ \in(R,2R); \   || D_{x}v||_{1} > M R^{-1-\frac{(a+2)q+(b+2)}{pq-1}}\},\\
\label{mDu}\Gamma_{4} (R)&:=&\{r\ \in(R,2R); \   ||D_{x} u||_{1}  > M R^{-1-\frac{(b+2)p+(a+2)}{pq-1}}\},\\
\label{mD2v}\Gamma_{5} (R)&:=&\{r\ \in(R,2R); \   ||D_{x}^{2} v||^{1+\epsilon}_{1+\epsilon}  > M R^{-2-\frac{(a+2)q+(b+2)}{pq-1}+b\epsilon}\},\\
\label{mD2u}\Gamma_{6}(R) &:=&\{r\ \in(R,2R); \   ||D_{x}^{2} u||^{1+\epsilon}_{1+\epsilon} > M R^{-2-\frac{(b+2)p+(a+2)}{pq-1} +a\epsilon}\}.
\end{eqnarray*}

Using (\ref{D2u}), we get 
\begin{eqnarray*}
 C &\ge&  R^{-N+2+\frac{(b+2)p+(a+2)}{pq-1}-a\epsilon} \int_{R}^{2R} ||D_{x}^{2} u||^{1+\epsilon}_{1+\epsilon} r^{N-1} dr \\
& \ge&  R^{-N+2+\frac{(b+2)p+(a+2)}{pq-1}-a\epsilon} |\Gamma_6(R)| R^{N-1} M  R^{-2-\frac{(b+2)p+(a+2)}{pq-1} +a\epsilon}= M |\Gamma_6(R)| R^{-1} .
 \end{eqnarray*} 
Therefore, choosing large enough $M$, we get $|\Gamma_{6}(R)|\le R/7$. Similarly, using (\ref{v})-(\ref{D2v}), one can see $|\Gamma_{i}(R)|\le R/7$ for $1\le i\le 5$. Hence, for each $R\ge 1$, we can find 
\begin{equation} \label{hatr}
 \hat R\in (R,2R)\setminus \bigcup_{i=1}^{i=6}\Gamma_{i}(R)\neq\phi.
\end{equation}
We now have the following upper bounds on (\ref{G1}) and (\ref{G2}) for the radius $\hat R$ given by (\ref{hatr});
\begin{eqnarray*}
G_{1}(\hat R)& \le&  C\ \hat R ^{N+a+2(p+1)} \left( \hat R^{  \left(-\frac{(a+2)q+(b+2)}{pq-1}-2+b\epsilon \right) \frac{1}{1+\epsilon} } + \hat R^{-2-\frac{(a+2)q+(b+2)}{pq-1} } \right)^{p+1}\\
&&  + C\ \hat R ^{N+b+2(q+1)} \left( \hat R^{  \left(-\frac{(b+2)p+(a+2)}{pq-1}-2+a\epsilon \right) \frac{1}{1+\epsilon} } + \hat R^{-2-\frac{(b+2)p+(a+2)}{pq-1} } \right)^{q+1},\\
&\le& C \left( \hat R^{-a_{1}(\epsilon)} +\hat R^{-a'_{1}(\epsilon)}   \right),
\end{eqnarray*}
where 
\begin{eqnarray*}
a_{1}(\epsilon)= (p+1) \left[  \left(2+\frac{(a+2)q+(b+2)}{pq-1}-b\epsilon\right) \frac{1}{1+\epsilon} -2- \frac{N+a}{p+1}\right],\\
a'_{1}(\epsilon)= (q+1) \left[  \left(2+\frac{(b+2)p+(a+2)}{pq-1}-a\epsilon\right) \frac{1}{1+\epsilon}  -2- \frac{N+b}{q+1}\right].\\
\end{eqnarray*}
Also, 
\begin{eqnarray*}
G_{2}(\hat R)& \le&  C\ \hat R ^{N+2} \left( \hat R^{  \left(-\frac{(b+2)p+(a+2)}{pq-1} -2+a\epsilon  \right) \frac{1}{1+\epsilon} } + \hat R^{-2-\frac{(b+2)p+(a+2)}{pq-1} } \right) \\
&&
 \left( \hat R^{  \left(-\frac{(a+2)q+(b+2)}{pq-1} -2+b\epsilon  \right) \frac{1}{1+\epsilon} } + \hat R^{-2-\frac{(a+2)q+(b+2)}{pq-1}} \right),
\\
&\le&  C \ \hat R^{-a_{2}(\epsilon)} ,
\end{eqnarray*}
where $$a_{2}(\epsilon)= -N-2+\frac{1}{1+\epsilon} \left(  4-(a+b)\epsilon  + \frac{(b+2)(p+1)+(a+2)(q+1)}{pq-1}  \right).$$

Hence, from (\ref{F}) we get
$$F( R)\le C \left(G_{1}(\hat R)+G_{2}(\hat R)\right)\le C \ R^{-\eta_\epsilon},$$ where $\eta_\epsilon:=\min\{ a_{1}(\epsilon), a'_{1}(\epsilon),a_{2}(\epsilon)\}$ and the positive constant $C$ does not depend on $R$. By a straightforward calculation, we have 
$$ a_{2}(0)= -N+2   + \frac{(b+2)(p+1)+(a+2)(q+1)}{pq-1}  >0 \ \ \text{iff} \ \ \frac{N+a}{p+1}+\frac{N+b}{q+1}>N-2.$$
Also, 
\begin {eqnarray}
\label{ineq1}a_{1}(0)>0, \ \ \text{iff} \ \ \frac{(a+2)q+(b+2)}{pq-1}>\frac{N+a}{p+1},\\
\label{ineq2}a'_{1}(0)>0, \ \ \text{iff} \ \ \frac{(b+2)p+(a+2)}{pq-1}>\frac{N+b}{q+1}.
\end{eqnarray}
Now, if $p$ and $q$ satisfy (\ref{curve}), then (\ref{ineq1}) and (\ref{ineq2}) hold, and we can therefore choose $\eta_\epsilon >0$ for small enough $\epsilon>0$. We now conclude by sending $R\to \infty$ and get the contradiction.

\hfill $\Box$ 

\section{On solutions of the second order H\'enon equation with finite Morse index}

We shall prove here Theorem \ref{resultMorse}. For that we  recall that a critical point $u\in C^{2} (\Omega)$ of the energy functional
\begin{equation*}
I(u):=\int_{\Omega} \frac{1}{2} |\nabla u|^{2}-\frac{1}{p+1} |x|^a u^{p+1}.
\end{equation*}
is said to be 
\begin{itemize}
 \item a {\it stable} solution of (\ref{mainsingle}) if  for any $\phi\in C_c^1(\Omega)$, we have
 $$I_{uu}(\phi):=\int_{\Omega} |\nabla \phi|^2 - p\int_{\Omega} |x|^a u^{p-1 }\phi^2\ge 0. $$

\item a  {\it stable solution outside a compact set} $\Sigma\subset\Omega$ if $I_{uu}(\phi)\ge 0$ for all $\phi\in C_c^1(\Omega\setminus\Sigma)$, also $u$ has a {\it Morse index} equal to $m\ge 1$ if $m$ is the maximal dimension of a subspace $X_m$ of $C^1_c(\Omega)$ such that $I_{uu}(\phi)<0$ for all $\phi\in X_m\setminus \{0\}$.

\item a {\it solution with Morse index $m$} if  there exist $\phi_1, ..., \phi_m$ such that $X_m= Span\{\phi_1, ..., \phi_m\}\subset C_c^1(\Omega)$ and $I_{uu}(\phi)<0$ for all $\phi\in X_m\setminus \{0\}$.
\end{itemize}

Note that if $u$ is of Morse index $m$, then  for all $\phi\in C_c^1(\Omega\setminus \Sigma)$ we have $I_{uu}(\phi)\ge 0$, where $\Sigma = \cup_{i=1}^{m} supp(\phi_i)$, and therefore $u$ is stable outside the  compact set $\Sigma\subset \Omega$. 
\\

We shall need the following lemma.
\begin{lemma}  \label{stpower}
Let $\Omega\subset \mathbb{R}^N$ and let $u\in C^2(\Omega)$ be a positive stable solution of (\ref{mainsingle}). Set $f(x)=|x|^a,a>0$, then, for any  $ 1\le t<-1+2p+2\sqrt{p(p-1)}$ we have
\begin{eqnarray} \label{estpower}
\int_{\Omega}\left( |\nabla u|^2 u^{t-1}+ f(x)  u^{t+p}\right)\phi^{2m} &\le&C \int_{\Omega} f(x)^{-\frac{t+1}{p-1}  }| \nabla \phi|^{2\frac{t+p}{p-1}},
\end{eqnarray}
 for all $\phi\in C_c^1(\Omega)$ with $0\le\phi\le 1$ and for large enough $m$. The constant $C$ does not depend on $\Omega$ and $u$.
\end{lemma}

{\bf Proof:} The following  proof also holds true for weak solutions. The ideas are adapted from \cite{egg,egg2,f}. Note first that for any stable solution of (\ref{mainsingle}) and $\eta\in C_c^1(\Omega)$, we have the following:
\begin{eqnarray}\label{stab}
p \int_\Omega  |x|^a u^{p-1} \eta^2 &\le&  \int_\Omega  |\nabla \eta|^2,\\
\label{pde} \int_\Omega |x|^a u^{p} \eta &=& \int_\Omega \nabla u\cdot\nabla \eta.
\end{eqnarray} 

Test (\ref{pde}) on $\eta=u^t \phi^2$ for $\phi\in C_c^1(\Omega)$ for an appropriate $t\in \mathbb{R}$ that will be chosen later,  to get 
\begin{eqnarray*}
\int_\Omega |x|^a u^{t+p} \phi^2 &=& \int_\Omega  \nabla u\cdot\nabla \left(u^t \phi^2\right)\\
& = &t \int_\Omega |\nabla u|^2 u^{t-1} \phi^2 + 2 \int_\Omega u^t \nabla u\cdot\nabla \phi \phi.
\end{eqnarray*} 
Apply Young's inequality\footnote{For any $a,b,\epsilon>0$, $ab\le \epsilon a^2+C(\epsilon) b^2$, for some $C(\epsilon)$.} to $\left(  |\nabla u| u^{\frac{t-1}{2}} \phi \right)\left( u^{\frac{t+1}{2}} |\nabla \phi|  \right)$ to obtain  
\begin{eqnarray}\label{pde10}
(t-\epsilon) \int_\Omega |\nabla u|^2 u^{t-1} \phi^2 \le C_{\epsilon} \int_\Omega  u^{t+1}|\nabla \phi|^2 + \int_\Omega |x|^a u^{t+p} \phi^2.
 \end{eqnarray} 
Now, test (\ref{stab}) on $u^{\frac{t+1}{2}}\phi$ to get
\begin{eqnarray*}
p \int_\Omega |x|^a u^{t+p} \phi^2 &\le& \frac{(t+1)^2}{4} \int_\Omega |\nabla u|^2  u^{t-1} \phi^2 + \int_\Omega u^{t+1}  |\nabla\phi|^2  \\
&+& (t+1) \int_\Omega  u^t  \nabla u \cdot \nabla \phi \phi \\
&\le& \left(   \frac{(t+1)^2}{4} +2\epsilon  \right) \int_\Omega |\nabla u|^2 u^{t-1} \phi^2 + (C'_{\epsilon,t}+C''_{\epsilon,t} ) \int_\Omega   u^{t+1}|\nabla \phi |^2,
 \end{eqnarray*} 
where again we have used Young's inequality in the last estimate.  Combine now this inequality with  (\ref{pde10}) to see
\begin{eqnarray} \label{majorest}
\left(   p -\frac{ \frac{(t+1)^2}{4} +2\epsilon  }{t-\epsilon}    \right)
\int_\Omega |x|^a u^{t+p} \phi^2 \le \left(   \frac{     \frac{(t+1)^2}{4} +2\epsilon  }{t-\epsilon} C_{\epsilon}  +C'_{\epsilon,t}+C''_{\epsilon,t}\right) \int_\Omega  u^{t+1}|\nabla \phi|^2.
 \end{eqnarray}
For an  appropriate choice of $t$, given in the assumption, we see that the coefficient in L.H.S. is positive for $\epsilon$ small enough. Therefore, replacing $\phi $ with $\phi^m$ for large enough $m$ and applying H\"{o}lder's  inequality  with exponents $\frac{t+p}{t+1}$ and $\frac{t+p}{p-1}$ we obtain 
\begin{eqnarray} \label{majorest1}
\int_\Omega |x|^a u^{t+p} \phi^{2m} \le D_{\epsilon,t,m} \int_\Omega |x|^{-\frac{t+1}{p-1}a}  |\nabla\phi |^{2\frac{t+p}{p-1}}. 
 \end{eqnarray}
Note that both exponents are greater than 1 for $t$ given in (i) and (ii).

On the other hand, combining (\ref{pde10}) and (\ref{majorest}) gives us
\begin{eqnarray*}
\int_\Omega |\nabla u|^2 u^{t-1} \phi^2 \le D'_{\epsilon,t} \int_\Omega  u^{t+1}|\nabla \phi|^2.
 \end{eqnarray*}
Similarly, replace $\phi $ by $\phi^m$ and apply H\"{o}lder's  inequality with exponents $\frac{t+p}{t+1}$ and $\frac{t+p}{p-1}$ to get  
\begin{eqnarray*}
\int_\Omega |\nabla u|^2 u^{t-1} \phi^{2m} \le D''_{\epsilon,t,m} \int_\Omega  |x|^{ -\frac{t+1}{p-1}a}   |\nabla \phi|^{ 2 \frac{t+p}{p-1}}.
 \end{eqnarray*}
This inequality and (\ref{majorest1}) finish the proof of (\ref{estpower}). 

\hfill $\Box$ 

Now, we are in the position to prove the theorem. 
\\
\\
{\bf Proof of Theorem \ref{resultMorse}:}  We proceed in the following  steps.
\\
\\
\textbf{Step 1}: We have the following standard Pohozaev type identity on any $\Omega\subset\mathbb{R}^N$. 
\begin{equation}\label{pohozaevidentity}
\frac{N+a}{p+1} \int_{\Omega} |x|^a u^{p+1} - \frac{N-2}{2} \int_{\Omega} |\nabla u|^2 = \frac{1}{p+1} \int_{\partial \Omega} |x|^a u^{p+1} x\cdot\nu + \int_{\partial\Omega} x\cdot\nabla u \nu\cdot\nabla u -\frac{1}{2} \int_{\partial\Omega} |\nabla u|^2 x\cdot\nu .
\end{equation}
 To get (\ref{pohozaevidentity}), just multiply both sides of (\ref{mainsingle}) by $x\cdot\nabla u$, do integration by parts and collect terms.
\\
\\
\textbf{Step 2}: The following estimates hold:
\begin{eqnarray*}
|\nabla u| &\in& L^2(\mathbb{R}^N),\\
|x|^a u^{p+1} &\in& L^1(\mathbb{R}^N).
\end{eqnarray*}

First recall that $u$ is stable outside a  compact set $\Sigma\subset \Omega$. To prove our claim, we use (\ref{estpower}) with the following test function $\xi_R\in C^1_c(\mathbb{R}^N\setminus \Sigma)$ for $R>R_0+3$ and $\Sigma \subset B_{R_0}$;
 $$\xi_R(x):=\left\{
                      \begin{array}{ll}
                      0, & \hbox{if $|x|<R_0+1$;} \\
                        1, & \hbox{if $R_0+2<|x|<R$;} \\
                        0, & \hbox{if $|x|>2R$;} 
                                                                       \end{array}
                    \right.$$
which satisfies $0\le\xi_R\le1$, $||\nabla\xi_R||_{L^{\infty}(B_{2R}\setminus B_{R})}<\frac{C}{R}$ and $||\nabla\xi_R||_{L^{\infty}(B_{R_0+2}\setminus B_{R_0+1})}<C_{R_0}$. Therefore, 
$$\int_{R_0+2<|x|<R} (|\nabla u|^2 u^{t-1} + |x|^a u^{t+p})\le C_{R_0}+\hat C\ R^{N- \frac{2(t+p)}{p-1}   -\frac{t+1}{p-1} a  }, $$ 
for all $1\le t<-1+2p+2\sqrt{p(p-1)}$. 
\\
Now,  set $t=1$ and send $R\to \infty$. Since $N< \frac{2(p+a+1)}{p-1}$, we see
$\int_{\mathbb{R}^N} |\nabla u|^2<\infty$ and $\int_{\mathbb{R}^N}|x|^a u^{p+1}<\infty$.
\\
\\
\textbf{Step 3}:  The following equality holds 
\begin{eqnarray}\label{equalitypoho}
\int_{\mathbb{R}^N} |\nabla u|^2 = \int_{\mathbb{R}^N}|x|^a u^{p+1}.
\end{eqnarray}

Multiply (\ref{mainsingle}) with $u\zeta_R$ for $\zeta_R\in C^1_c(\mathbb{R}^N)$ which satisfies  $0\le\zeta_R\le1$, $||\nabla\zeta_R||_{\infty}<\frac{C}{R}$ and
 $$\zeta_R(x):=\left\{
                      \begin{array}{ll}
                        1, & \hbox{if $|x|<R$;} \\
                        0, & \hbox{if $|x|>2R$.} 
                                                                       \end{array}
                    \right.$$
Then, integrate over $B_{2R}$ to get 
\begin{eqnarray} \label{identitytest}
\int_{B_{2R}} |x|^a u^{p+1} \zeta_R - \int_{B_{2R}} |\nabla u|^2 \zeta_R = \int_{B_{2R}} \nabla \zeta_R\cdot \nabla u \ u.
\end{eqnarray}
By H\"{o}lder's inequality, we have the following upper bound for R.H.S. of (\ref{identitytest}),
\begin{eqnarray*}
|\int_{B_{2R}}  \nabla \zeta_R\cdot \nabla u \ u| & \le & R^{-1} \int_{B_{2R}} |\nabla u| ( |x|^{\frac{a}{p+1}}u )\ |x|^{-\frac{a}{p+1}}\\
&\le& R^{-1}\left(\int_{B_{2R}}  |\nabla u|^2 \right)^{\frac{1}{2}}      \left(\int_{B_{2R}}  |x|^a u^{p+1} \right)^{\frac{1}{p+1}}  \left(\int_{B_{2R}}  |x|^{-\frac{2a }{p-1}} \right)^{\frac{p-1}{2(p+1)}}\\
&=& R^{\frac{N(p-1)}{2(p+1)}-\frac{a }{p+1}-1}\left(\int_{B_{2R}}  |\nabla u|^2 \right)^{\frac{1}{2}}      \left(\int_{B_{2R}}  |x|^a u^{p+1} \right)^{\frac{1}{p+1}}.
\end{eqnarray*}
Therefore, from Step 2, there exists a positive constant $C$ independent of $R$ such that 
\begin{eqnarray*}
|\int_{B_{2R}}  \nabla \zeta_R\cdot \nabla u \ u| & \le & C \ R^{\frac{N(p-1)-2(a+p+1)}{2(p+1)}}.
\end{eqnarray*}
Since $N< \frac{2(p+a+1)}{p-1}$, we have $\lim_{R\to \infty }|\int_{B_{2R}} \nabla \zeta_R\cdot \nabla u \ u| =0$. Hence (\ref{identitytest}) implies (\ref{equalitypoho}).
\\
\\
\noindent\textbf{Step 4}: we have $$(\frac{N+a}{p+1}-\frac{N-2}{2})\int_{\mathbb{R}^N} |x|^a u^{p+1}=0.$$  

Apply Lemma \ref{stpower} for $t=1$ with the following test function $\phi_R\in C^1_c(\mathbb{R}^N\setminus \Sigma)$ for $R>2 R_0$;

 $$\phi_R(x):=\left\{
                      \begin{array}{ll}
                      0, & \hbox{if $|x|<R/2$;} \\
                        1, & \hbox{if $R<|x|<2R$;} \\
                        0, & \hbox{if $|x|>3R$;} 
                                                                       \end{array}
                    \right.$$
which satisfies $0\le\phi_R\le1$, $||\nabla\phi_R||_{L^{\infty}(B_{3R}\setminus B_{R/2})}<\frac{C}{R}$ to get 

\begin{eqnarray} \label{stpowerineq}
\int_{B_{2R}\setminus B_{R}}\left( |\nabla u|^2 + |x|^a  u^{p+1}\right) &\le&C R^{N-\frac{2(p+a+1)}{p-1}}.
\end{eqnarray}
Now, define the following sets for large enough $M$;
\begin{eqnarray*}
\label{mv}\theta_{1}(R) &:=&\{r\ \in(R,2R); \  || D_x u(r)||^2_{2}> M R^{-\frac{2(p+a+1)}{p-1}}\},\\
\label{mu}\theta_{2} (R)&:=&\{r\ \in(R,2R); \  || u(r)||^{p+1}_{p+1}  > M R^{-\frac{2(p+a+1)}{p-1}-a} \}.
\end{eqnarray*}
From (\ref{stpowerineq}), we have
\begin{eqnarray*}
 C &\ge&  R^{-N+\frac{2(p+a+1)}{p-1}+a} \int_{R}^{2R}|| u(r)||^{p+1}_{p+1} r^{N-1} dr \\
& \ge&  R^{-N+\frac{2(p+a+1)}{p-1}+a} |\theta_2(R)| R^{N-1} M  R^{-\frac{2(p+a+1)}{p-1}-a}= M |\theta_2(R)| R^{-1} .
 \end{eqnarray*} 
Similarly, one can show $|\theta_1(R)|\le R/M$. By choosing $M$ large enough we conclude $|\theta _i(R)| \le {R}/{3}$ for $i=1,2$. Therefore, for each $R\ge 1$, we can find 
\begin{equation*}
 \tilde R\in (R,2R)\setminus \bigcup_{i=1}^{i=2}\Lambda_{i}(R)\neq\phi.
\end{equation*}

Now, apply Pohozaev identity, (\ref{pohozaevidentity}), with $\Omega=B_{\tilde R}$ to see that R.H.S. converges to zero if $R\to \infty$ for subcritical $p$, i.e. $N<\frac{2(p+a+1)}{p-1}$. Hence, 
\begin{eqnarray*} 
\frac{N-2}{2}\int_{\mathbb{R}^N} |\nabla u|^2 = \frac{N+a}{p+1}\int_{\mathbb{R}^N} |x|^a u^{p+1}.
\end{eqnarray*}

From this and (\ref{equalitypoho}), we finish the proof of Step 4.

\hfill $\Box$

 \noindent \textbf{Remark}: For the Sobolev critical case $p=\frac{N+2+2a}{N-2}$, using the change of variable $w:=u(r^{1+\frac{a}{2}})$ and applying well-known classifying-type results mentioned in the introduction for the Lane-Emden equation, one can see all radial solutions of (\ref{mainsingle}) are of the following form 
\begin{equation}\label{radial}
u_{\epsilon}(r):= k(\epsilon) (\epsilon+r^{2+a})^{\frac{2-N}{2+a}},
\end{equation}
where $k(\epsilon)=\left(\epsilon(N+a)({N-2})\right)^{\frac{N-2}{2(2+a)}}$. Then, from the classical Hardy's inequality it is  straightforward to see $u_\epsilon$ is stable outside a compact set $\overline{B_{R_0}}$, for an appropriate $R_0$. Note that for $-2<a\le 0$, by Schwarz symmetrization (or rearrangement), it is shown in \cite{gy} that all radial solutions of (\ref{mainsingle}) with $p=\frac{N+2+2a}{N-2}$ and $N>2$ are of the form (\ref{radial}).

\section{On solutions of the fourth order H\'enon equation with finite Morse index}

We shall prove here Theorem \ref{fourthresultMorse}. For that we  recall that a critical point $u$ of the energy functional
\begin{equation*}
I(u):=\int_{\Omega} \frac{1}{2} |\Delta u|^{2}-\frac{1}{p+1}\int_{\Omega}  |x|^a u^{p+1},
\end{equation*}
is said to be a {\it stable} solution of (\ref{fourth}),  if  for any $\phi\in C_c^4(\Omega)$, we have
 $$I_{uu}(\phi):=\int_{\Omega} |\Delta \phi|^2 - p\int_{\Omega} |x|^a u^{p-1 }\phi^2\ge 0. $$
Similarly to the second order case, one can define the notion of {\it stability outside a compact set}, which contains the notion of solutions with {\it finite Morse index}. We first prove the following estimate.

\begin{lemma}  \label{fourthstpower}
Let $\Omega\subset \mathbb{R}^N$ and let $u\in C^4(\Omega)$ be a positive stable solution of (\ref{fourth}). Then, for large enough $m$, we have for all $\phi\in C_c^4(\Omega)$ with $0\le\phi\le 1$, 
\begin{eqnarray} \label{estpowerfourth}
\int_{\Omega}\left( |\Delta u|^2 + |x|^a  u^{p+1}\right)\phi^{2m} &\le&C \int_{\Omega} |x|^{-\frac{2}{p-1}  a}| T (\phi)|^{\frac{p+1}{p-1}},
\end{eqnarray}
 where $T(\phi):=|\Delta \phi|^2+|\nabla\phi|^4+|\Delta|\nabla\phi|^2|+|\nabla \phi\cdot\nabla\Delta\phi|$. The constant $C$ does not depend on $\Omega$ and $u$.
\end{lemma}

\textbf{Proof:} For any stable solution of (\ref{fourth}) and $\eta\in C_c^4(\Omega)$, we have the followings:
\begin{eqnarray}\label{fourthstab}
p \int_\Omega  |x|^a u^{p-1} \eta^2 &\le&  \int_\Omega  |\Delta \eta|^2,\\
\label{fourthpde} \int_\Omega |x|^a u^{p} \eta &=& \int_\Omega \Delta u\Delta \eta.
\end{eqnarray} 

Test (\ref{fourthpde}) on $\eta=u \phi^2$ for $\phi\in C_c^4(\Omega)$  to get 
\begin{eqnarray} \label{fourthst}
\int_\Omega |x|^a u^{p+1} \phi^2 &=& \int_\Omega  \Delta u\ \Delta \left(u \phi^2\right)
\end{eqnarray} 

Also, test (\ref{fourthstab}) on $u \phi$ and use (\ref{fourthst}) to get
\begin{eqnarray*}
(p-1) \int_\Omega |x|^a u^{p+1} \phi^2 &\le&  \int_\Omega |\Delta( u\phi)|^2 - \int_\Omega |x|^a u^{p+1} \phi^2 \\
&=&  \int_\Omega |\Delta( u\phi)|^2- \int_\Omega \Delta u\Delta (u\phi^2).
 \end{eqnarray*} 

By a straightforward calculation, one can see that the following identity holds:
\begin{eqnarray}\label{identity}
|\Delta (u\phi)|^2 - \Delta u\Delta (u\phi^2) = 4 |\nabla u\cdot\nabla \phi|^2 + u^2 |\Delta \phi|^2-2u\Delta u |\nabla \phi|^2 +2 \nabla u^2\cdot\nabla \phi\Delta \phi.
\end{eqnarray} 
Therefore, we have
\begin{eqnarray*}
(p-1) \int_\Omega |x|^a u^{p+1} \phi^2 &\le&  4\int_\Omega |\nabla u|^2|\nabla \phi|^2 + \int_\Omega u^2 |\Delta \phi|^2-2\int_\Omega u\Delta u |\nabla \phi|^2 \\ 
&&+ 2\int_\Omega\nabla u^2\cdot\nabla \phi\Delta \phi.
 \end{eqnarray*} 
A simple integration by parts yields
\begin{eqnarray}\label{integration}
\int_\Omega |\nabla u|^2|\nabla \phi|^2 = \int_\Omega  u(-\Delta u) |\nabla \phi|^2   + \frac{1}{2}\int_\Omega u^2 \Delta |\nabla \phi|^2,
 \end{eqnarray} 
which then simplifies the previous inequality to become 
\begin{eqnarray*}
(p-1) \int_\Omega |x|^a u^{p+1} \phi^2 &\le&  6 \int_\Omega u(-\Delta u) |\nabla \phi|^2 + \int_\Omega u^2( -|\Delta \phi|^2+2\Delta |\nabla\phi|^2-2\nabla\phi\cdot\nabla \Delta \phi).
 \end{eqnarray*} 
Therefore,
\begin{eqnarray}\label{majorineq}
 \int_\Omega |x|^a u^{p+1} \phi^2 &\le&  C \int_\Omega u|\Delta u| |\nabla \phi|^2 + \int_\Omega u^2L(\phi),
 \end{eqnarray} 
where $L(\phi):= |\Delta \phi|^2+2|\Delta |\nabla\phi|^2|+2|\nabla\phi\cdot\nabla \Delta \phi|$.

On the other hand, from (\ref{identity}) and (\ref{integration}), one can see
\begin{eqnarray*}
\int_\Omega |\Delta (u\phi)|^2& =&\int_\Omega \Delta u\Delta (u\phi^2) + 4\int_\Omega |\nabla u\cdot\nabla \phi|^2 +\int_\Omega u^2 |\Delta \phi|^2-2u\Delta u |\nabla \phi|^2 -2\int_\Omega  u^2div(\nabla \phi\Delta \phi)\\
&=& \int_\Omega |x|^a u^{p+1} \phi^2 + 6 \int_\Omega u(-\Delta u) |\nabla \phi|^2 +\int_\Omega u^2 \left(-|\Delta \phi|^2+2\Delta |\nabla\phi|^2-2\nabla\phi\cdot\nabla \Delta \phi\right).
 \end{eqnarray*} 
By combining (\ref{majorineq}), the identity $\Delta (u\phi)=\phi\Delta u+2\nabla u\cdot\nabla \phi+u\Delta\phi $ and Young's inequality,  we get the following estimate 
\begin{eqnarray*}
\int_\Omega |\Delta u|^2\phi^2\le C \int_\Omega u|\Delta u| |\nabla \phi|^2 + C \int_\Omega u^2L(\phi).
 \end{eqnarray*} 
Therefore, 
\begin{eqnarray*}
 \int_\Omega \left(|x|^a u^{p+1} +|\Delta u|^2\right)\phi^{2} \le C  \int_\Omega u|\Delta u| |\nabla \phi|^2+C \int_\Omega u^2 L(\phi).
  \end{eqnarray*} 
Now, replacing $\phi$ with $\phi ^m$ for large enough $m>0$ and applying Young's inequality we end up with 
\begin{eqnarray*}
 \int_\Omega \left(|x|^a u^{p+1} +|\Delta u|^2\right)\phi^{2m}& \le& C  \int_\Omega u|\Delta u| |\nabla \phi|^2\phi^{2(m-1)}+ C\int_\Omega u^2 L(\phi^m)\\
 &\le& \epsilon  \int_\Omega |\Delta u|^2 \phi^{2m} +C_\epsilon \int_\Omega u^2 |\nabla \phi|^4 \phi^{2(m-2)}+C \int_\Omega u^2 L(\phi^m).
  \end{eqnarray*}
Then, for large enough $m$ 
\begin{eqnarray}\label{finalineq}
 \int_\Omega \left(|x|^a u^{p+1} +|\Delta u|^2\right)\phi^{2m}& \le&  C \int_\Omega u^2\phi^{2(m-2)} T(\phi) ,
  \end{eqnarray}
where $T(\phi):=|\Delta \phi|^2+|\nabla\phi|^4+|\Delta|\nabla\phi|^2|+|\nabla \phi\cdot\nabla\Delta\phi|$. Now, apply H\"{o}lder's inequality to get
\begin{eqnarray*}
\int_\Omega u^2 \phi^{2(m-2)} T(\phi) &=&  \int_\Omega |x|^\frac{2a}{p+1}u^2 \phi^{2(m-2)} |x|^{-\frac{2a}{p+1}} T(\phi)  \\
&\le& \left( \int_\Omega |x|^a u^{p+1}  \phi^{2(m-2)\frac{p+1}{2}} \right)^\frac{2}{p+1} \left(    \int_\Omega  |x|^{-\frac{2a}{p-1}} T^\frac{p+1}{p-1}(\phi)  \right)^\frac{p-1}{p+1}
 \end{eqnarray*} 
 Choosing $m$ large enough, say $2(m-2)\frac{p+1}{2}\ge 2m$, from (\ref{finalineq}) we finally get the desired inequality
 \begin{eqnarray*}
 \int_\Omega (|x|^a u^{p+1} +|\Delta u|^2)\phi^{2m} \le C \int_\Omega |x|^{-\frac{2a}{p-1}}  T^\frac{p+1}{p-1}(\phi).
  \end{eqnarray*}

\hfill $\Box$ 
\\
\\
{\bf Proof of Theorem \ref{fourthresultMorse}:} We proceed in the following  steps.
\\
\\
\textbf{Step 1}: We have the following standard Pohozaev type identity on any $\Omega\subset\mathbb{R}^N$. 
\begin{eqnarray}\label{fourthpohozaevidentity}
\nonumber\frac{N+a}{p+1} \int_{\Omega} |x|^a u^{p+1} - \frac{N-4}{2} \int_{\Omega} |\Delta u|^2 &=& \frac{1}{p+1} \int_{\partial \Omega} |x|^a u^{p+1} x\cdot\nu  -\frac{1}{2} \int_{\partial\Omega} |\Delta u|^2 x\cdot\nu \\
&&-\int_{\partial\Omega} \nabla\Delta u\cdot\nu x\cdot\nabla u + \int_{\partial\Omega} \Delta u \nabla(x\cdot\nabla u)\cdot \nu.
\end{eqnarray}
 To get (\ref{fourthpohozaevidentity}), just multiply both sides of (\ref{fourth}) by $x\cdot\nabla u$, do integration by parts and collect terms.
\\
\\
\textbf{Step 2}: we have
\begin{eqnarray*}
|\Delta u| &\in& L^2(\mathbb{R}^N),\\
|x|^a u^{p+1} &\in& L^1(\mathbb{R}^N).
\end{eqnarray*}

Since $u$ is stable outside a compact set $\Sigma\subset \Omega$, using (\ref{estpowerfourth}) with the following test function $\xi_R\in C^1_c(\mathbb{R}^N\setminus \Sigma)$ for $R>R_0+3$ and $\Sigma \subset B_{R_0}$;
 $$\xi_R(x):=\left\{
                      \begin{array}{ll}
                      0, & \hbox{if $|x|<R_0+1$;} \\
                        1, & \hbox{if $R_0+2<|x|<R$;} \\
                        0, & \hbox{if $|x|>2R$;} 
                                                                       \end{array}
                    \right.$$
which satisfies $0\le\xi_R\le1$, $||D^i\xi_R||_{L^{\infty}(B_{2R}\setminus B_{R})}<\frac{C}{R^i}$ and $||D^i\xi_R||_{L^{\infty}(B_{R_0+2}\setminus B_{R_0+1})}<C_{R_0}$ for $i=1,\cdots,4$, we get 
$$\int_{R_0+2<|x|<R} (|\Delta u|^2  + |x|^a u^{p+1})\le C_{R_0}+\hat C\ R^{N- \frac{4(p+1)}{p-1}   -\frac{2}{p-1} a  }. $$ 

 For subcritical exponents, $N< \frac{2(2p+a+2)}{p-1}$, we see
$\int_{\mathbb{R}^N} |\Delta u|^2<\infty$ and $\int_{\mathbb{R}^N}|x|^a u^{p+1}<\infty$.
\\
\\
\textbf{Step 3}:  The following equality holds 
\begin{equation}\label{equality}
\int_{\mathbb{R}^N} |x|^a u^{p+1}=\int_{\mathbb{R}^N} |\Delta u|^2.
\end{equation}

Multiply (\ref{fourth}) with $u\zeta_R$ for $\zeta_R\in C^4_c(B_{2R})$ which satisfies  $0\le\zeta_R\le1$, $||D^i\zeta_R||_{\infty}<\frac{C}{R^i}$ for $i=1,\dots,4$ and
 $$\zeta_R(x):=\left\{
                      \begin{array}{ll}
                        1, & \hbox{if $|x|<R$;} \\
                        0, & \hbox{if $|x|>2R$.} 
                                                                       \end{array}
                    \right.$$
Then, integrate over $B_{2R}$ to get 
\begin{eqnarray} \label{fourthidentitytest}
\int_{B_{2R}} |x|^a u^{p+1} \zeta_R - \int_{B_{2R}} |\Delta u|^2 \zeta_R =\int_{B_{2R}} u \Delta u\Delta\zeta_R  + 2 \int_{B_{2R}} \Delta u\nabla u\cdot \nabla \zeta_R=:I_1(R)+I_2(R) .
\end{eqnarray}
By H\"{o}lder's inequality, we have the following upper bound for $I_1(R)$,
\begin{eqnarray*}
|I_1(R)| & \le & R^{-2} \int_{B_{2R}} |\Delta u| ( |x|^{\frac{a}{p+1}}u )\ |x|^{-\frac{a}{p+1}}\\
&\le& R^{-2}\left(\int_{B_{2R}}  |\Delta u|^2 \right)^{\frac{1}{2}}      \left(\int_{B_{2R}}  |x|^a u^{p+1} \right)^{\frac{1}{p+1}}  \left(\int_{B_{2R}}  |x|^{-\frac{2a }{p-1}} \right)^{\frac{p-1}{2(p+1)}}\\
&=& R^{\frac{N(p-1)}{2(p+1)}-\frac{a }{p+1}-2}\left(\int_{B_{2R}}  |\Delta u|^2 \right)^{\frac{1}{2}}      \left(\int_{B_{2R}}  |x|^a u^{p+1} \right)^{\frac{1}{p+1}}.
\end{eqnarray*}
Therefore, from Step 2, there exists a positive constant $C$ independent of $R$ such that 
\begin{eqnarray*}
|I_1(R)| & \le & C \ R^{\frac{N(p-1)-2(a+p+1)}{2(p+1)}}.
\end{eqnarray*}
Since $N< \frac{2(p+a+1)}{p-1}$, we have $\lim_{R\to \infty }|I_1(R)| =0$. Now, we consider the second term in R.H.S. of (\ref{fourthidentitytest}). Apply Young's inequality for a given $\epsilon>0$ (we choose it later) to get 
$$|I_2(R)|\le \epsilon \int_{\mathbb{R}^N}  |\Delta u|^2 + C_\epsilon \int_{B_{2R}} |\nabla u|^2|\nabla \zeta_R|^2,$$
Using Green's theorem we get 
$$\int_{B_{2R}} |\nabla u|^2|\nabla \zeta_R|^2=\int_{B_{2R}}u(-\Delta u) |\nabla\zeta_R|^2+\frac{1}{2} \int_{B_{2R}} u^2 \Delta |\nabla \zeta_R|^2=:I_3(R)+I_4(R).$$

By the same discussion as given for $I_1(R)$ one can see $\lim_{R\to \infty }|I_3(R)| =0$. For the term $I_4(R)$, we apply H\"{o}lder's inequality again
\begin{eqnarray*}
|I_4| & \le & R^{-4} \int_{B_{2R}} |x|^{\frac{a}{p+1}}u^2 \ |x|^{-\frac{a}{p+1}}\\
&\le& R^{-4}      \left(\int_{B_{2R}}  |x|^a u^{p+1} \right)^{\frac{2}{p+1}}  \left(\int_{B_{2R}}  |x|^{-\frac{2a }{p-1}} \right)^{\frac{p-1}{p+1}}\\
&=& R^{\frac{N(p-1)}{(p+1)}-\frac{2a }{p+1}-4}    \left(\int_{B_{2R}}  |x|^a u^{p+1} \right)^{\frac{2}{p+1}}.
\end{eqnarray*}
By Step 2 and sending $R$ to infinity we get, $\lim_{R\to \infty }|I_4(R)| =0$. Since $\lim_{R\to \infty }|I_2(R)|\le \epsilon \int_{\mathbb{R}^N} |\Delta u|^2$ for any $\epsilon >0$, we have $\lim_{R\to \infty }|I_2(R)|=0$. Therefore, (\ref{equality}) follows.
\\
\\
\textbf{Step 4}:  The following equality holds $$(\frac{N+a}{p+1}-\frac{N-4}{2})\int_{\mathbb{R}^N} |x|^a u^{p+1}=0.$$

Apply Lemma \ref{fourthstpower} with the following test function $\phi_R\in C^1_c(\mathbb{R}^N\setminus \Sigma)$ for $R>2 R_0$;

 $$\phi_R(x):=\left\{
                      \begin{array}{ll}
                      0, & \hbox{if $|x|<R/2$;} \\
                        1, & \hbox{if $R<|x|<2R$;} \\
                        0, & \hbox{if $|x|>3R$;} 
                                                                       \end{array}
                    \right.$$
where $0\le\phi_R\le1$, $||D^i\phi_R||_{L^{\infty}(B_{3R}\setminus B_{R/2})}<\frac{C}{R^i}$. Then, we get 

\begin{eqnarray} \label{fourthestpower}
\int_{B_{2R}\setminus B_{R}} |\Delta u|^2 + |x|^a  u^{p+1} &\le&C R^{N-\frac{2(2p+2+a)}{p-1}}.
\end{eqnarray}
On the other hand, we are interested in similar upper bounds for the following terms 
$$J_1(R):=\int_{B_{2R}\setminus B_{R}} |\Delta u| |\nabla u|  \ \ \ \ \text{and } \ \ \ \ \ J_2(R):=\int_{B_{2R}\setminus B_{R}} |\Delta u| |D^2_x u| .$$
For the first term, $J_1(R)$, using Schwarz's inequality we have
\begin{equation*}
\int_{B_{2R}\setminus B_{R}} |\Delta u| |\nabla u| < \left( \int_{B_{2R}\setminus B_{R}} |\Delta u|^2 \right)^{1/2}  \left( \int_{B_{2R}\setminus B_{R}} |\nabla u|^2     \right)^{1/2}.
\end{equation*}
From standard elliptic interpolation estimates, L$^2$-norm version of Lemma \ref{interp},  we have 
\begin{eqnarray*}
\int_{B_{2R}\setminus B_{R}}  |  \nabla u |^2 &\le& C R^2  \int_{B_{4R}\setminus B_{R/2}}  |\Delta u|^2  +  C  R^{-2}     \int_{B_{4R}\setminus B_{R/2}}    u^2  \\
&\le& C R^{N-\frac{2(2p+2+a)}{p-1}+2}+ R^{\frac{N(p-1)}{(p+1)}-\frac{2a }{p+1}-2}    \left(\int_{\mathbb{R}^N}  |x|^a u^{p+1} \right)^{\frac{2}{p+1}}\\
&= & C R^{\frac{p-1}{p+1}\left(N-\frac{2(2p+2+a)}{p-1}+2\right)+2} \left( R^{\frac{2}{p+1}\left(N-\frac{2(2p+2+a)}{p-1}\right) } +\left(\int_{\mathbb{R}^N}  |x|^a u^{p+1} \right)^{\frac{2}{p+1}} \right)
\end{eqnarray*}
Since $\int_{\mathbb{R}^N}  |x|^a u^{p+1}<\infty$ and $N<\frac{2(2p+2+a)}{p-1}$, for $R>1$ we have
\begin{eqnarray*}
\int_{B_{2R}\setminus B_{R}}  |  \nabla u |^2 \le C R^{\frac{p-1}{p+1}\left(N-\frac{2(2p+2+a)}{p-1}+2\right)+2}
\end{eqnarray*}
Therefore, 
\begin{equation}\label{fourthnabla}
\int_{B_{2R}\setminus B_{R}} |\Delta u| |\nabla u| <C R^{\frac{p}{p+1}\left(N-\frac{2(2p+2+a)}{p-1}\right)+1}.
\end{equation}
Similarly for the second term, $J_2(R)$, using Lemma \ref{ellip}, i.e.,
$$\int_{B_{2R}\setminus B_{R}}    |  D_{x}^2 u |^2\le C\left(\int_{B_{4R}\setminus B_{R/2}   }  |\Delta u|^2  +   R^{-4}   \int_{B_{4R}\setminus B_{R/2}   }    u^2  \right ),$$
and similar type discussions one can see
\begin{equation}\label{fourthsecond}
\int_{B_{2R}\setminus B_{R}}  |\Delta u| |D^2_x u|  <C  R^{\frac{p}{p+1}\left(N-\frac{2(2p+2+a)}{p-1}\right)}.
\end{equation}

Now, define the following sets for large enough $M$;
\begin{eqnarray*}
\label{mv}\Lambda_{1}(R) &:=&\{r\ \in(R,2R); \  || \Delta_x u(r)||^2_{2}> M R^{-\frac{2(2p+2+a)}{p-1}}\},\\
\label{mu}\Lambda_{2} (R)&:=&\{r\ \in(R,2R); \  || u(r)||^{p+1}_{p+1}  > M R^{-\frac{2(2p+2+a)}{p-1}-a} \},\\
\label{mDv}\Lambda_{3} (R)&:=&\{r\ \in(R,2R); \   || \Delta_x u(r) |\nabla_x u(r)|||_{1} > M R^{-\frac{p}{p+1}\left(\frac{N}{p}+\frac{2(2p+2+a)}{p-1}\right)+1}\},\\
\label{mDu}\Lambda_{4} (R)&:=&\{r\ \in(R,2R); \   ||\Delta_x u(r) D^2_{x}u(r) ||_{1}  > M R^{-\frac{p}{p+1}\left(\frac{N}{p}+\frac{2(2p+2+a)}{p-1}\right)}\}.
\end{eqnarray*}
In the following, we shall find a bound for the measure of the above sets. From (\ref{fourthestpower}), we have
\begin{eqnarray*}
 C &\ge&  R^{-N+\frac{2(2p+2+a)}{p-1}+a} \int_{R}^{2R}|| u(r)||^{p+1}_{p+1} r^{N-1} dr \\
& \ge&  R^{-N+\frac{2(2p+2+a)}{p-1}+a} |\Lambda_2(R)| R^{N-1} M  R^{-\frac{2(2p+2+a)}{p-1}-a}= M |\Lambda_2(R)| R^{-1} .
 \end{eqnarray*} 
Also, from (\ref{fourthnabla}) 
\begin{eqnarray*}
 C &\ge& R^{\frac{p}{p+1}\left(-N+\frac{2(2p+2+a)}{p-1}\right)-1}\int_{R}^{2R} || \Delta_x u(r) |\nabla_x u(r)|||_{1}  r^{N-1} dr \\
& \ge&  R^{\frac{p}{p+1}\left(-N+\frac{2(2p+2+a)}{p-1}\right)-1} |\Lambda_3(R)| R^{N-1} M  R^{-\frac{p}{p+1}\left(\frac{N}{p}+\frac{2(2p+2+a)}{p-1}\right)+1}= M |\Lambda_3(R)| R^{-1} .
 \end{eqnarray*} 
Similarly, from (\ref{fourthsecond}) and (\ref{fourthestpower}) we get $|\Lambda_1(R)|,|\Lambda_4(R)|\le R/M$. By choosing $M$ large enough we conclude $|\Lambda _i(R)| \le {R}/{5}$ for $i=1,\cdots,4$. Therefore, for each $R\ge 1$, we can find 
\begin{equation} \label{tilder}
 \tilde R\in (R,2R)\setminus \bigcup_{i=1}^{i=4}\Lambda_{i}(R)\neq\phi.
\end{equation}
Then, from the definition of $\tilde R$ and $\Lambda_i$ for $i=1,\cdots,4$, we have 
\begin{eqnarray}\label{ineq1}
\int_{|x|=\tilde R}  |\Delta_x u(\tilde R)| |D^2_x u (\tilde R)   |  &\le&C  \tilde {R}^{\frac{p}{p+1}\left(N-\frac{2(2p+2+a)}{p-1}\right)-1}\\ \label{ineq2}
\int_{|x|=\tilde R}  |\Delta_x u  (\tilde R) | |\nabla_x u  (\tilde R)  |  &\le&C  \tilde {R}^{\frac{p}{p+1}\left(N-\frac{2(2p+2+a)}{p-1}\right)}\\ \label{ineq3}
\int_{|x|=\tilde R}  |\Delta_x u (\tilde R)   |^2 &\le&C  \tilde {R}^{\frac{p}{p+1}\left(N-\frac{2(2p+2+a)}{p-1}\right)-1}\\   \label{ineq4}
\int_{|x|=\tilde R} u^{p+1}(\tilde R)  &\le&C  \tilde {R}^{\frac{p}{p+1}\left(N-\frac{2(2p+2+a)}{p-1}\right)-a-1}
 \end{eqnarray} 

Using (\ref{fourthpohozaevidentity}) with $\Omega=B_{2\tilde R}\setminus B_{\tilde R}$, one can see
\begin{eqnarray}\label{ineq5}\left| \int_{|x|=\tilde R} \nabla \Delta u\cdot\nu x\cdot\nabla u  \right|< C \tilde{R}^{\frac{p}{p+1}\left(N-\frac{2(2p+2+a)}{p-1}\right)}.  \end{eqnarray}

Now, applying the Pohozaev identity, (\ref{fourthpohozaevidentity}), with $\Omega= B_{\tilde R}$ and using (\ref{ineq1})-(\ref{ineq5}),  R.H.S. of (\ref{fourthpohozaevidentity}), converges to zero if $R\to \infty$ for subcritical $p$, i.e. $N<\frac{2(2p+2+a)}{p-1}$. Hence, 
\begin{eqnarray*} \label{identitypohoz}
\frac{N-4}{2}\int_{\mathbb{R}^N} |\Delta u|^2 = \frac{N+a}{p+1}\int_{\mathbb{R}^N} |x|^a u^{p+1}.
\end{eqnarray*}
From this and (\ref{equality}), we finish the proof of Step 4.

\hfill $\Box$

\section{On stable solutions of the H\'enon-Lane-Emden system}

We shall prove here Theorem \ref{system}. For that we recall that a classical solution $(u,v)$ of (\ref{mainbound}) is said to be {\it pointwise stable} if there exists positive smooth $\zeta,\eta$ such that 
\begin{eqnarray}
\label{pointstable}
 \left\{ \begin{array}{lcl}
\hfill -\Delta \zeta&=& p|x|^{a}v^{p-1}\eta   \ \ \text{in}\ \ \Omega,\\   
\hfill -\Delta \eta&=& q|x|^{b}u^{q-1}\zeta   \ \ \text{in}\ \ \Omega.
\end{array}\right.
  \end{eqnarray}
In what follows we give the stability inequality for system (\ref{mainbound}). This inequality is the novelty here and is key tool in proving Theorem 4. The idea of geting such an inequality comes from \cite{fg}. 

\begin{lemma} \label{sysineq}
Assume that $(u,v)$ is a pointwise stable solution of (\ref{mainbound}), then for any test function $\phi\in C_c^1(\Omega)$, we have
\begin{equation}\label{sysstable}
\sqrt{pq}\int_{\Omega} |x|^{\frac {a+b}{2}  } v^{\frac{p-1}{2}}  u^{  \frac{q-1}{2}} \phi^2 \le \int_{\Omega} |\nabla \phi|^2.
\end{equation}
\end{lemma}

\noindent{\bf Proof:}  Let $(u,v)$ be a pointwise stable solution of (\ref{mainbound}) in such a way that there exists positive smooth $\zeta,\eta$ such that (\ref{pointstable}).  Multiply the first equation by $\phi^2 \zeta^{-1}$ and the second equation by  $\phi^2 \eta^{-1}$, integrate by parts and use Young's inequality to get 
\begin{eqnarray*}
p \int_{B_R} |x|^a v^{p-1} \frac{\eta}{\zeta}\phi^2&=&-\int_{B_R} \frac{\Delta \zeta}{\zeta}\phi^2\le  \int_{B_R} |\nabla \phi|^2\\
q \int_{B_R} |x|^b u^{q-1} \frac{\zeta}{\eta}\phi^2&=&-\int_{B_R} \frac{\Delta \eta}{\eta}\phi^2\le \int_{B_R} |\nabla \phi|^2.
\end{eqnarray*}
Adding these two equations and doing simple calculations we get 
\begin{eqnarray*}
2 \int_{B_R} |\nabla \phi|^2 \ge \int_{B_R} \left(p|x|^a v^{p-1} \frac{\eta}{\zeta}+q |x|^b u^{q-1} \frac{\zeta}{\eta}    \right)\phi^2\ge 2\sqrt{pq}\int_{B_R} |x|^{\frac{a+b}{2}} u^{\frac{q-1}{2}} v^{\frac{p-1}{2}} \phi^2.
\end{eqnarray*}
\hfill $\Box$ 

The following pointwise estimate is taken from \cite{phan}. As was said before, the first version of this paper was done independently of \cite{phan} and without using the following lemma of Phan. This last section --which uses Lemma \ref{phancompare}--  was added after his paper was posted. 

\begin{lemma}\label{phancompare}[Phan, \cite{phan}] Assume that $(u,v)$ is a classical  solution for  (\ref{mainbound}), then for 
\begin{equation}\label{ab}
 0\le a-b\le (N-2)(p-q)
 \end{equation} we have 
\begin{equation}\label{compare}
|x|^av^{p+1}  \le \frac{p+1}{q+1} |x|^b u^{q+1}.
 \end{equation}
\end{lemma}

Combining the above lemmas we conclude the following integral estimate which is a counterpart of Lemma \ref{stpower} for the second order case and Lemma \ref{fourthstpower} for the fourth order case.

\begin{lemma}  \label{sysstpower}
For $\Omega\subset \mathbb{R}^N$, assume that (\ref{ab}) holds and that $(u,v)$ is a pointwise stable solution of (\ref{mainbound}). Set $\theta:=\frac{pq(q+1)}{p+1}$. Then, for any $t$ such that 
 $$\sqrt{\theta}-\sqrt{\theta-\sqrt{\theta}}<t<\sqrt{\theta}+\sqrt{\theta-\sqrt{\theta}},$$ we have for all $\phi\in C_c^2(\Omega)$ such that $0\le\phi\le 1$,
\begin{eqnarray} \label{sysestpower}
\int_{\Omega} |x|^a  v^p u^{2t-1} \phi^{2} &\le&C \int_{\Omega} u^{2t} \left(| \nabla \phi|^{2}+|\Delta \phi|\right).
\end{eqnarray}
The constant $C$ does not depend on $\Omega$ and $(u,v)$.
\end{lemma}

\noindent{\bf Proof:} Note first that for $p\ge q$, we have $\theta\ge q^2>1$ and also
$$\frac{1}{2}<  \sqrt{\theta}-\sqrt{\theta-\sqrt{\theta}}<1<\sqrt{\theta}+\sqrt{\theta-\sqrt{\theta}}.$$
 Let $(u,v)$ is a pointwise stable solution of (\ref{mainbound}). Then, Lemma \ref{sysineq} applies and by replacing $\phi$ with $u^{t}\phi$  in (\ref{sysstable}),  where $\phi$ is a test function, we obtain
\begin{equation}\label{sysstable1}
\sqrt{pq}\int_{} |x|^{\frac {a+b}{2}  } v^{\frac{p-1}{2}}  u^{  \frac{q-1}{2}} u^{2t} \phi^2 \le \int_{} |\nabla (u^t\phi)|^2.
\end{equation}
Rewriting the left hand side as
 $\sqrt{pq} \int |x|^{\frac {a+b}{2}  } v^{\frac{p-1}{2}}  u^{  \frac{q+1}{2}} u^{2t-1} \phi^2 $ 
 and using Lemma \ref{phancompare}, i.e. $\sqrt\frac{q+1}{p+1}|x|^{\frac{a-b}{2}} v^{\frac{p+1}{2}}\le u^{\frac{q+1}{2}}$, we get 
\begin{equation}\label{sysinteg}
\sqrt{\frac{pq(q+1)}{p+1}} \int |x|^a   v^p  u^{2t-1} \phi^2 \le t^2 \int |\nabla u|^2 u^{2t-2}\phi^2+\int  u^{2t}\phi |\Delta \phi|.
\end{equation}
To find an upper bound for the first term in  the above inequality with the gradient term, we multiply both sides of the first equation in (\ref{mainbound}) to get
\begin{eqnarray*}
\int_{} |x|^a  v^p u^{2t-1} \phi^{2}&=&\int \nabla u\cdot \nabla (u^{2t-1}\phi^2)\\&=&
(2t-1)\int |\nabla u|^2 u^{2t-2} \phi^2 -\frac{1}{2t} \int u^{2t} \Delta (\phi^2).
\end{eqnarray*}
Since $t> \frac{1}{2}$, we get 
$$t^2 \int |\nabla u|^2 u^{2t-2} \phi^2 \le \frac{t^2}{2t-1} \int_{} |x|^a  v^p u^{2t-1} \phi^{2}+C_t\int u^{2t}\left (\phi|\Delta \phi|+|\nabla \phi|^2\right).  $$
Combining this and (\ref{sysinteg}) we have 
\begin{equation}\label{sysinteg}
\left(\sqrt{\frac{pq(q+1)}{p+1}} -  \frac{t^2}{2t-1}  \right)\int |x|^a   v^p  u^{2t-1} \phi^2 \le C_t\int u^{2t}\left (\phi|\Delta \phi|+|\nabla \phi|^2\right).
\end{equation}
\hfill $\Box$

\noindent{\bf Proof of Theorem \ref{system}:} Define $z:=u^\tau$ for $1<2\sqrt{\theta}-2\sqrt{\theta-\sqrt{\theta}}<\tau<2\sqrt{\theta}+2\sqrt{\theta-\sqrt{\theta}}$ and $\theta:=\frac{pq(q+1)}{p+1}$. Then, 
\begin{equation*}
|\Delta z|\le  C \left( |\nabla u|^2 u^{\tau -2} + |x|^a u^{\tau -1} v^p \right).
\end{equation*} 
By integrating over balls we get 
\begin{equation}\label{deltaz}
\int_{B_R} |\Delta z|\le  C \int_{B_R} |\nabla u|^2 u^{\tau -2} + C\int_{B_R} |x|^a u^{\tau -1} v^p .
\end{equation} 

We are now after an upper bound for the right hand side of the above inequality. To control the second term, apply Lemma \ref{stpower} for  $t:=\frac{\tau}{2}$ and standard test function $\zeta_R$ used in the proof of Theorem \ref{fourthresultMorse} to get 
\begin{eqnarray*} 
\int_{B_R} |x|^a  v^p u^{\tau-1} &\le&C R^{-2} \int_{B_R} u^{\tau}.
\end{eqnarray*}
To bound the first term, we use the first equation of the system. Multiply both sides of (\ref{mainbound}) with $u^{\tau-1}\zeta_R^2$ and integrate by parts to get 

\begin{eqnarray*} 
\int_{B_R}  |\nabla u|^2 u^{\tau -2} & \le& \int_{B_R}  |\nabla u|^2 u^{\tau -2} \zeta_R^2
\\&=& \frac{1}{\tau-1} \int_{B_R} |x|^a  v^p u^{\tau-1} \zeta_R^2+ \frac{1}{\tau(\tau-1)} \int_{B_R} u^\tau \Delta(\zeta_R^2)
\\&\le& C\int_{B_R} |x|^a  v^p u^{\tau-1}+C R^{-2} \int_{B_R} u^\tau.
\end{eqnarray*}
Therefore, the following upper bound holds for (\ref{deltaz}),
\begin{equation*}
\int_{B_R} |\Delta z|\le C R^{-2} \int_{B_R} u^\tau,
\end{equation*} 
which means $ R^2 || \Delta z||_{L^1(B_{R})} \le C ||z||_{L^1(B_{R})} $. Now, applying Lemma \ref{regularity}  for $z=u^\tau$ we get 
$$ || z||_{L^k(B_R)}\le C R^{N (\frac{1}{k}-1) } ||z||_{L^1(B_{2R})},$$
where $C=C(k,N)>0$ and any $1\le k<\frac{N}{N-2}$.

Now take $1\le k_i<\frac{N}{N-2}$ for $1\le i\le n$ and $2\sqrt{\theta}-2\sqrt{\theta-\sqrt{\theta}}<2t:=\tau k_{n-1}!<2\sqrt{\theta}+2\sqrt{\theta-\sqrt{\theta}}$. The notation "!" stands for $k_{n-1}!:= \prod_{i=0}^{n-1} k_i$ and set $k_0=1$. By induction  we have
$$ || z||_{L^{k_n!}(B_R)}\le C R^{\tilde k_n } ||z||_{L^1(B_{2R})},$$
where $\tilde k_n=N\sum_{i=1}^{n} \frac{1-k_i}{k_i!}=N\left( \frac{1}{k_n!}-1 \right)$ and $C=C(k_i,N)>0$. So, 
$$ \left(\int_{B_R} u^{\tau k_n!} \right)^{\frac{1}{k_n!}}\le C R^{ N\left( \frac{1}{k_n!}-1 \right) }\int_{B_{2R}} u^{\tau}  .$$

Let $0<\tau<q$ and from Corollary \ref{L1est}  we get 
$$ \int_{B_{2R}} u^\tau \le C R^{N-\frac{p(b+2)+a+2}{pq-1}\tau}.$$
Therefore
\begin{eqnarray} \label{tau}
 \left(\int_{B_R} u^{\tau k_n!} \right)^{\frac{1}{k_n!}}&\le& C R^{\tau \left(  \frac{N}{\tau k_n!}- \frac{p(b+2)+a+2}{pq-1}  \right)}.
\end{eqnarray}
 So, in the following dimensions 
\begin{eqnarray*} 
  N <\frac{p(b+2)+a+2}{pq-1}  \tau k_n! 
  \end{eqnarray*}
the right hand side of (\ref{tau}) converges to zero as $R$ tends to infinity. Note that since $\tau k_{n-1}!<2\sqrt{\theta}+2\sqrt{\theta-\sqrt{\theta}}$ and $k_n<\frac{N}{N-2}$, we have $\tau k_n!<(2\sqrt{\theta}+2\sqrt{\theta-\sqrt{\theta}})\frac{N}{N-2}$. So,
\begin{eqnarray*} 
  N <\tau k_n! \frac{p(b+2)+a+2}{pq-1}  <  2+  \frac{p(b+2)+a+2}{pq-1} \left(2\sqrt{\theta}+2\sqrt{\theta-\sqrt{\theta}}\right)
  \end{eqnarray*}
  Recall that $\theta:=\frac{pq(q+1)}{p+1}$, which completes the proof.

\hfill $\Box$

\end{document}